\documentclass[12pt]{article}
\usepackage{graphicx}
\usepackage[final]{pdfpages}
\usepackage{amsmath}
\usepackage{amsfonts}
\usepackage{amssymb}
\usepackage{subfigure}
\usepackage{algorithm,algorithmicx,algpseudocode}
\usepackage{pgfplots}
\usepackage{placeins}
\usepackage{caption}
\usepackage{colonequals}
\usepackage{microtype}

\providecommand{\U}[1]{\protect\rule{.1in}{.1in}}
\providecommand{\U}[1]{\protect\rule{.1in}{.1in}}
\providecommand{\U}[1]{\protect\rule{.1in}{.1in}}
\newtheorem{theorem}{Theorem}

\numberwithin{equation}{section} \numberwithin{theorem}{section}
\newtheorem{definition}[theorem]{Definition}
\newtheorem{example}[theorem]{Example}

\newcommand{\Addresses}{{
  \medskip
  \footnotesize

  \noindent$^1$ Optimization Department, Fraunhofer - ITWM, Kaiserslautern 67663, Germany\\
  $^2$ Mathematics Department, ORT Braude College, Karmiel
2161002, Israel\par\nopagebreak

\medskip

  \noindent E-mail: esther.bonacker@itwm.fraunhofer.de and avivg@braude.ac.il
}}

\begin{document}

\title{Speedup of lexicographic optimization by superiorization and its
applications to cancer radiotherapy treatment}
\author{Esther Bonacker$^1$, Aviv Gibali$^2$, Karl-Heinz K\"{u}fer$^1$ \\
and Philipp S\"{u}ss$^1$}
\date{September 12, 2016}
\maketitle \Addresses

\begin{abstract}
Multicriteria optimization problems occur in many real life applications,
for example in cancer radiotherapy treatment and in particular in intensity modulated
radiation therapy (IMRT). In this work we
focus on optimization problems with multiple objectives that are ranked
according to their importance. We solve these problems numerically by combining lexicographic optimization
with our recently proposed level set scheme, which yields a sequence of
auxiliary convex feasibility problems; solved here via projection methods.
The projection enables us to combine the newly introduced superiorization methodology
with multicriteria optimization methods to speed up computation while guaranteeing convergence of the optimization.
We demonstrate our scheme with a simple 2D academic example (used in the literature) and also present results from calculations on four real head neck cases in IMRT (Radiation Oncology of the Ludwig-Maximilians University, Munich, Germany) for two different choices of superiorization parameter sets suited to yield fast convergence for each case individually or robust behavior for all four cases.

\end{abstract}



\section{Introduction}

A wide variety of real-life problems involve the simultaneous
optimization of several objectives or criteria. Objectives that contradict
each are common, thus resulting in situations where one
objective cannot be improved without worsening at least one other
objective. For example, in cancer radiotherapy treatment, one goal is to
maximize the amount of cancerous tissue exposed to radiation treatment while
another (contradicting) goal is to minimize the radiation exposure of healthy
organs. Here, there is typically no trivial compromise.
Problems of this kind are called \textit{Multi-Criteria Optimization} (MCO)
problems.\vspace\baselineskip

The number of competing goals in the treatment planning for \textit{intensity modulated radiation therapy} (IMRT) is quite large (about 10-25). Exploring the solutions to the MCO problem involves a high-dimensional Pareto boundary \cite{geoffrion68}. Approximating this set is computationally expensive, and experience shows that there are many tradeoffs that are either clinically unacceptable or irrelevant. Thus a significant part of the computational effort is wasted on identifying useless solutions.

The complexity of the approximated high-dimensional Pareto boundary (see e.g. \cite{kmssaabt08}) makes it difficult for the decision maker to navigate. An alternative approach would be to first identify the most important aspects of the IMRT plan and make sure that these are met by the solution; compromises, if necessary, would only be acceptable for subsidiary goals.

Toward this end, we want to incorporate a priori knowledge about priorities into the treatment plan optimization process. This information translates naturally into \textit{Lexicographic Optimization} (LO) problems, which constitute a special class of MCO problems.\vspace\baselineskip

Classical LO methods minimize the objective functions
sequentially, starting with the most important one and proceeding according
to the lexicographic order of the objectives. Here, the optimal
value found for each objective is added as a constraint for subsequent
optimizations.

These constraints can be quite rigorous, however, depending on the correlation of the corresponding functions and the decision maker might not always want to sacrifice so much optimizing freedom for the sake of a single goal. We therefore propose the more intuitive approach of gathering goals in groups of descending importance, as in \cite{jmf07}, and formulating objective functions as weighted sums of the functions contained in these groups. This reduces the number of objective functions to be minimized, which leads to fewer optimization problems to be solved and thus lowers the computational effort.\vspace\baselineskip

Typically, the subsequent objectives are in conflict with each
other, which implies that the minimum of one objective function is a
relatively bad starting point for subsequent optimization levels.

To address this issue, we propose to combine the superiorization methodology with an MCO algorithm.
Superiorization can speed up convergence by steering the
MCO algorithm towards solutions that are better suited to minimizing
subsequent objective functions. In the best case, one could even find the
minimizer of all subsequent objective functions before reaching the last
optimization level. Here, subsequent minimization becomes unnecessary, and one need only verify that the solution is
indeed the minimizer of all subsequent objective functions under the given constraints.\vspace\baselineskip

The superiorization methodology \cite{cdh10, pscr10, hgdc12, jcj13,
dcsgx15} was recently developed as a framework for algorithms that lie
conceptually between feasibility-seeking and optimization algorithms. It is
a heuristic tool that does not guarantee to find the optimum value of a
given functional, rather it obtains a solution that is superior (with
respect to a given objective function) to the solution achieved by a
classical feasibility seeking algorithmic operator.

The main advantage of an algorithm that uses superiorization as the driving
tool (as opposed to optimization) is that it requires less computational
resources - given a good choice of parameters - while providing comparable
solutions, from the point of view of real-world applications, to those that
one would get with algorithms that use optimization.

Within this framework, we wish to include the class of projection methods,
which are applied successfully in many real-world feasibility problems, see e.g., \cite{bk13}.
This class has witnessed great progress in recent
years and its member algorithms have been applied with success to problems in sensor networks, image reconstruction, image
processing, IMRT planning, resolution enhancement and in many others; see \cite{pscr10, cccdhh12, dcsgx15}. Apart
from theoretical interest, the main advantage of projection methods is their computational efficiency. They
commonly have the ability to handle huge-size problems of dimensions in
which other, more sophisticated methods cease to be
efficient.

In this paper we propose to apply the superiorization methodology with
projection methods to improve the speed of convergence of lexicographic optimization
methods. We demonstrate the plausibility of using superiorization with a simple
example with two optimization variables. We also present results from calculations on four real head neck cases in IMRT for two different choices of superiorization
parameter sets suited to yield fast convergence for each case individually
or robust behavior for all four cases.

The paper is organized as follows. In Section \ref{Sec:Methods} we present
some preliminaries, definitions and algorithms that will be needed in the
sequel. Later, in Section \ref{Sec:Numerical illustration} our general
scheme is illustrated on a 2D example and we present the results of our
calculations on the aforementioned IMRT head neck cases. Finally we discuss
our findings and conclusions in Sections \ref{Sec:Discussion} and
\ref{Sec:Conclusion}.

\section{Methods}

\label{Sec:Methods}

\subsection{Problems, algorithms and theory}

Let us start by giving a definition of what is a multicriteria
optimization problem.

\begin{definition}
Given a mapping $F=(f_{1},\ldots ,f_{m}):\mathbb{R}^{n}\rightarrow \mathbb{R}%
^{m}$ with $f_{i}:\mathbb{R}^{n}\rightarrow \mathbb{R}$ for all $i\in
I=\{1,\ldots ,m\}$ and let $\left\{ g_{_{j}}:\mathbb{R}^{n}\rightarrow
\mathbb{R}\right\} _{j\in J}$ for $J=\{1,\ldots ,s\}$. The \texttt{
Multicriteria Optimization} (MCO) problem is the following.

\begin{align}
\text{\textquotedblleft min\textquotedblright }\hspace{0.1cm}& F(x)  \notag
\\
\text{such that }x& \in \bigcap_{j\in J}\Omega _{j}  \label{problem:MCO}
\end{align}
where for all $j\in J$, $\Omega _{j}\colonequals\left\{ x\in \mathbb{R}^{n}\mid
g_{_{j}}(x)\leq 0\right\} $.
\end{definition}

\subsubsection{Lexicographic optimization}

In this work we focus on a particular MCO
method which is called \textit{lexicographic optimization }(LO). For LO we
define an order for the objectives according to their importance and
minimize them in that order under the given constraints. After having
minimized one objective, we impose an additional constraint which ensures
that the objective function value of subsequent solutions do not deviate
from the previously found minimum by more than a user-defined value.

In the context of IMRT we consider a more general setting for the
formulation of the optimization problems. Here, we gather the objective
functions $f_{1},\ldots ,f_{m}$ in $M$ priority groups which are ranked from
most to least important (compare e.g., \cite{jmf07}). Denote by $I_{\mu }\subseteq I$ an ordered subset
of objective function indices such that $\{I_{\mu }\}_{\mu =1}^{M}$ is a
\textit{partition} of $I$, meaning that

\begin{equation}
\bigcup_{\mu =1}^{M}I_{\mu }=I\text{ and } I_{\mu }\cap I_{\gamma }=\emptyset \quad \text{for } 1\leq\mu\neq \gamma\leq M.
\end{equation}

It is clear that if for all $\mu =1,\ldots,M$ the set $I_{\mu }$ is a singleton then this reformulation reduces to the classical lexicographic optimization
problem.

Now we define%
\begin{equation}\label{eq:Phi_mu}
\phi _{\mu }(x)\colonequals\sum_{i\in I_{\mu }}w_{i}f_{i}(x)
\end{equation}
where for each $i\in I_{\mu }$, $w_{i}$ is the weight of
the $i$-th objective function $f_{i}$. So $\phi_{\mu }$ gathers the objective functions whose indices are elements of $I_{\mu }$. At optimization level $\mu$ $(\leq M)$ we thus solve the following problem.

\begin{align}\label{eq:MCLO}
 \min \text{ } &\phi _{\mu }(x)  \notag \\
\text{such that } & x\in \bigcap_{j\in J}\Omega _{j} \notag\\
    & \phi _{\gamma }(x)\leq \phi _{\gamma }^{\ast }+\delta _{\gamma }\quad
\gamma =1,...,\mu -1
\end{align}
where $\phi _{\gamma }^{\ast }$ are the minimum values found for each $\phi
_{\gamma }$ during the previous optimization levels, $\delta
_{\gamma }\geq 0$ are small user-chosen constants for $\gamma =1,...,\mu -1$ and $\Omega _{j}$ are as in (\ref{problem:MCO}).

Let us illustrate the above for a two-stage LO also known as \textit{%
bi-level optimization}. Assume that $I_{1}=\{1\},$ $I_{2}=\{2\}$, and $%
w_{1}=w_{2}=1$; Thus $\phi _{1}=f_{1}$ and $\phi _{2}=f_{2}$. At the first
optimization level we solve the following single-objective optimization
problem.%
\begin{align}
\min \text{ } & f_{1}(x)  \notag \\
\text{such that }& x \in \bigcap_{j\in J}\Omega _{j}.  \label{eq:p1}
\end{align}%
After we have solved (\ref{eq:p1}) and obtained $f_{1}^{\ast }$ ($=\phi _{1}^{\ast }$), we continue
with the second optimization level, where our problem formulates as
\begin{align}\label{eq:p2}
 \min \text{ } & f_{2}(x)  \notag \\
\text{such that }& x\in \bigcap_{j\in J}\Omega _{j} \notag\\
    & x\in \Omega ^{\delta ,1}
\end{align}
where $\Omega ^{\delta ,1}\colonequals\left\{ x\in \mathbb{R}^{n}\mid
f_{1}(x)\leq f_{1}^{\ast }+\delta _{1}\right\} $ for some small
user-chosen constant $\delta _{1}\geq 0$.\vspace\baselineskip

In the above, the choice of $\delta_1$ and in general $\delta_{\mu}$ for $\mu=1,..., M$, defines in a fixed way how far from previously found optima of objective functions of higher priority the algorithm is allowed to deviate in the current optimization level. A different approach has been taken by \cite{lmfftr12} which allows the decision maker to choose the tradeoff between two subsequent (in terms of their lexicographic order) objective functions in an interactive way.\vspace\baselineskip

Observe that if the constraint set of (\ref{eq:p1})

\begin{equation}
\bigcap_{j\in J}\Omega _{j}\label{eq:bilevel}
\end{equation}
(or of (\ref{eq:p2}) with $\delta_1=0$) contains only one element, then this is the optimal solution of the bi-level
optimization problem and it is found already at the first optimization
level. Otherwise, in case of multiple feasible points, a feasibility seeking
algorithm for finding an element of (\ref{eq:p1}) might end up with a point which is
relatively far from the optimal solution of the second optimization level,
this is due to the fact that the second objective function has not been
taken into account yet. This effect is likely to occur when the objective
functions of the subsequent optimization levels represent conflicting goals
as commonly happens in cancer therapy.

Motivated by these difficulties we propose to steer the algorithm towards
solutions that provide good starting points for subsequent optimization
levels by applying the superiorization methodology together with projection
methods. Before we can do so we need to transform the optimization problem such that we can use projection methods to solve it.

\subsubsection{Level set scheme}

For the aforementioned purpose we use the \textit{level set scheme} \cite{gks} which is a new projection-based scheme for solving convex optimization problems.
It transforms a convex optimization problem into a sequence of auxiliary feasibility problems by iteratively constraining the objective function from above until the
feasibility problem is inconsistent. To illustrate the idea let $\phi:\mathbb{R}^{n}\rightarrow \mathbb{R}$ be convex and continuously
differentiable function and $\Xi\subseteq\mathbb{R}^{n}$ nonempty, closed and convex set. Following the level set scheme we can
reformulate the problem
\begin{align}
\min\text{ } & \phi (x)  \notag \\
\text{ such that }& x\in \Xi
\end{align}
in the following form

\begin{align}\label{eq:EpiForm}
 \min \text{ } &t  \notag \\
\text{such that }& \phi (x)\leq t \notag\\
& x\in \Xi.
\end{align}

Now let $\left\{  \varepsilon_{k}\right\}  _{k=0}^{\infty}$ be some user chosen positive sequence, for
example $\varepsilon_{k}\equiv\varepsilon=0.1\ $or $\varepsilon_{k}%
=0.1\left\vert \phi(x^{k})\right\vert $.

In the initial step of the level set scheme we choose any projection method to
find a feasible point $x^{0}\in \Xi $, set $t_{0}=\phi(x^{0})-\varepsilon_{0}$ and $k=0$.
Now at the iterative step $k>0$, when we are given the current point $x^{k-1}$, we try to solve the following \textit{convex feasibility problem}
(CFP)

\begin{align}
\text{Find a point }& x^{k}\in\mathbb{R}^n\notag\\
\text{ such that } & \phi (x^{k})\leq t_{k} \notag\\
& x^{k}\in \Xi.
\end{align}

If there exists a feasible solution, set $t_{k+1}=\phi(x^{k})-\varepsilon_{k}$
and continue, else there exists no feasible solution, and $x^{k}$ is an $\varepsilon_{k}%
$-optimal solution.

\subsubsection{Projection methods}

In this work, we use \textit{projection methods} to solve the auxiliary convex feasibility problems generated by the level set scheme. Projection methods are iterative algorithms that use
projections onto sets. They rely on the general principle that when a
family of (usually closed and convex) sets is present, projections onto
the given individual sets are easier to perform than projections onto other
sets (intersections, image sets under some transformation, etc.) that are
derived from the given individual sets.
Projection methods come in various algorithmic structures, some of which are
particularly suitable for parallel computing, and they demonstrate nice
convergence properties, for example bounded perturbations resilience. This
fact allows us to use them to solve the CFPs resulting from the level set
scheme and to incorporate superiorization, which is presented in the following paragraph, while retaining convergence to the optimal point.

Consider the CFP:
\begin{equation}
\text{Find a point }x\in \bigcap_{j\in J}\Omega _{j}
\end{equation}
with $\Omega _{j}$ defined as in \eqref{problem:MCO}. The particular method we use in this paper is the \textit{simultaneous subgradient projections method} which is defined as follows. Let $x^0\in\mathbb{R}^n$ be an arbitrary starting point.
Given the current iterate $x^{k}$, calculate the next iterate $x^{k+1}$ via

\begin{equation}
x^{k+1}\colonequals x^{k}-\lambda _{k}\sum_{j\in J,\text{ }g_{j}(x^{k})>0}w_{j}\frac{%
g_{j}(x^{k})}{\left\Vert \xi ^{k}\right\Vert ^{2}}\xi ^{k}
\end{equation}
where $\xi ^{k}\in \partial g_{j}(x^{k})$ (subgradient of $g_{j}$ at $x^{k}$
), $w_{j}>0$ are \textit{weights} and $\lambda _{k}\in \lbrack \epsilon _{1},2-\epsilon _{2}]$ (relaxation parameters) for arbitrary $\epsilon _{1},\epsilon _{2}>0$.

Simultaneous projection methods are also referred to as \textquotedblleft
parallel\textquotedblright\ methods. In this case, in order to evaluate the
next iterate, all (or a block of more than one) constraints are taken into
account. In case of a system of linear equalities, this method is known as
Cimmino's method \cite{Cimmino38}. See \cite{im86,
dos_Santos87} for more details.

\begin{definition}
Given a problem $\mathbb{P}$, an algorithmic operator $
\boldsymbol{T}:\mathbb{R}^{n}\rightarrow\mathbb{R}^{n}$ is said to be
\texttt{bounded perturbations resilient} if the following is true.\vspace\baselineskip

If the sequence $\{x^{k}\}_{k=0}^{\infty},$ generated by $x^{k+1}=%
\boldsymbol{T}(x^{k})$, for all $k\geq0$, converges to a solution of $\mathbb{P}$, then any sequence $\{y^{k}\}_{k=0}^{\infty}$ of points in $%
\mathbb{R}^{n}$ generated by $y^{k+1}=\boldsymbol{T}(y^{k}+\beta_{k}v^{k}),$
for all $k\geq0$, also converges to a solution of $\mathbb{P}$ provided
that, for all $k\geq0$, $\beta_{k}v^{k}$ are bounded perturbations,
meaning that $\beta_{k}\geq0$ for all $k\geq0$ such that ${\displaystyle\sum
\limits_{k=0}^{\infty}}\beta_{k}\,<\infty$ and the sequence $\{v^{k}\}_{k=0}^{\infty}$ is bounded.
\end{definition}

\subsubsection{Superiorization methodology}

The \textit{superiorization methodology} of \cite{hgdc12} was recently developed as a framework for algorithms that lie conceptually between feasibility-seeking and optimization
algorithms. It is designed to find a solution to a CFP which is superior
with respect to a given objective function $\psi $, meaning with a value of $%
\psi $ at least as low, but possibly lower, compared to the
solution obtained by a classical feasibility seeking algorithmic
operator. The state of current research on superiorization can
best be appreciated from the \textquotedblleft Superiorization and
Perturbation Resilience of Algorithms: A Bibliography compiled and
continuously updated by Yair Censor\textquotedblright\ which is
at: http://math.haifa.ac.il/yair/bib-superiorization-censor.html.
In particular, \cite{herman-review-sm} and \cite{weak-strong15}
are recent reviews of interest.

\begin{definition}
Given a function $\psi :\mathbb{R}^{n}\rightarrow \mathbb{R}$ and a point $
z\in \mathbb{R}^{n}$, we say that a vector $d\in \mathbb{R}^{n}$ is \texttt{nonascending} for $\psi $ at $z$ if and only if $\left\Vert d\right\Vert
\leq 1$ and there is a $\zeta >0$ such that for all $\lambda \in \lbrack
0,\zeta ]$ we have $\psi (z+\lambda d)\leq \psi (z).$
\end{definition}

Consider again the convex feasibility problem:

\begin{equation}
\text{Find a point }x\in \Xi  \label{eq:CCO}
\end{equation}
where $\Xi \subseteq \mathbb{R}^{n}$ is a closed and convex set and let $%
\boldsymbol{T}$ be any feasibility-seeking algorithmic operator, that
defines the basic algorithm $x^{k+1}=\boldsymbol{T}(x^{k})$ to solve the
convex feasibility problem.

Let $\psi :\mathbb{R}^{n}\rightarrow \mathbb{R}$ be a given convex and
continuously differentiable function. The \textit{superiorized} (with
respect to $\psi $) version of $\boldsymbol{T}$ is

\begin{equation}
y^{k+1}=\boldsymbol{T}(y^{k}+\beta _{k}d^{k})
\end{equation}
where $d^{k}$ is nonascending for $\psi $ at $y^{k}$.

One option for $d^{k}$, as it is used in \cite{hgdc12}, is
\begin{equation}
d^{k}=\left\{
\begin{array}{cc}
-\frac{\displaystyle\nabla \psi (y^{k})}{\displaystyle\left\Vert \nabla \psi
(y^{k})\right\Vert } & \text{if }\nabla \psi (y^{k})\neq 0 \\
0 & \nabla \psi (y^{k})=0.%
\end{array}
\right.\label{eq:nabla}
\end{equation}

If the solution set of the feasibility problem contains several elements,
the superiorized version of $\boldsymbol{T}$ obtains a solution that is
superior with respect to $\psi $ to the solution achieved by the basic
algorithmic operator $\boldsymbol{T}$.

\subsubsection{Main result}

We now combine all of the mentioned concepts to derive our approach. We consider the multicriteria optimization problem (\ref{problem:MCO}) and
define the priority groups $I_{\mu }$ and weights $\{w_{i}\}_{i\in I_{\mu
}}$, $\mu =1,...,M$. We
then transform it into a LO problem like (\ref{eq:MCLO}). Now we apply the
level set scheme and consequently solve the following optimization problem
at every optimization level $\mu $:

\begin{align}\label{CFPLexOpt}
 \min \text{ } &t^{(\mu )}  \notag \\
\text{such that }& x\in \Omega ^{t^{(\mu )}} \notag\\
    & x\in {\textstyle\bigcap\limits_{j\in J}}\Omega _{j} \notag\\
    & x\in {\textstyle\bigcap\limits_{\gamma =1}^{\mu -1}}\Omega ^{\delta
,\gamma }
\end{align}
where

\begin{align}
\Omega ^{t^{(\mu )}}&\colonequals\left\{ x\in \mathbb{R}^{n}\mid \phi_{\mu}(x)\leq
t^{(\mu )}\right\} \smallskip \notag \\
\Omega _{j}&\colonequals\left\{ x\in \mathbb{R}^{n}\mid g_{j}(x)\leq
0\right\}  \notag\\
\Omega ^{\delta,\gamma }&\colonequals\left\{ x\in
\mathbb{R}^{n}\mid \phi_{\gamma}(x)\leq \phi_{\gamma}^{\ast }(x)+\delta
_{\gamma }\right\}. \notag
\end{align}

To do that, we transform (\ref{CFPLexOpt}) into a sequence of auxiliary CFPs
of the form:

\begin{align}\label{eq:MCLOauxCFP}
\text{Find a point }& x^{k}\in\mathbb{R}^{n}\notag\\
\text{ such that } & x^{k}\in \Omega ^{t_{k}^{(\mu )}} \notag\\
    & x^{k}\in {\textstyle\bigcap\limits_{j\in J}}\Omega _{j} \notag\\
    & x^{k}\in {\textstyle\bigcap\limits_{\gamma =1}^{\mu -1}}\Omega ^{\delta
,\gamma }
\end{align}
for decreasing sequences of $\{t_{k}^{(\mu )}\}$. We solve the auxiliary
CFPs by using the simultaneous subgradient projections method as the basic
algorithm $\boldsymbol{T}$.\vspace\baselineskip

Let $\mathcal{S}\subseteq \textstyle%
\bigcup_{\gamma = \mu+1}^{M} I_{\gamma}$ be a subset of indices contained in
the subsequent priority groups $I_{\mu +1},\cdots ,I_{M}$ and let $\hat{w}%
_{i}$ be weights to the corresponding objective functions.

Then, after having successfully solved $K$ auxiliary CFPs we include
superiorization with respect to $\psi (x)\colonequals\sum_{i\in \mathcal{S}}\hat{w}%
_{i}f_{i}$ by defining $d^{k}$ as in (\ref{eq:nabla}) with $x^k$ and setting

\begin{equation}
x^{k+1}=x^{k}+\beta _{k}d^{k}
\end{equation}%
where $\beta _{k}\geq 0$ for all $k\geq 0$ such that ${\displaystyle%
\sum\limits_{k=0}^{\infty }}\beta _{k}\,<\infty $ and $\psi (x^{k}+\beta
_{k}d^{k})\leq \psi (x^{k})$, meaning that $d^{k}$ is a direction of
nonascend for $\psi $ at $x^{k}$.\vspace\baselineskip

In the following we refer to this method as \textit{superiorized lexicographic optimization} (SLO).

\subsection{Application to IMRT}

The IMRT planning problem is to find a treatment plan, i.e.\ energy
fluence intensities, which results in a dose distribution in the patient's
body that irradiates the tumor as homogeneously as possible while sparing
critical healthy organs.

The dose (measured in Gy) delivered to the patient's body is evaluated according to
clinical goals which are translated into the functions given below. The
clinical goals are ordered according to their importance and can also be
grouped together which reduces the number of optimization levels and
therefore accelerates the optimization process. For our calculations on four
head neck cases we chose one group of goals as constraints and three
subsequent priority groups.

\subsubsection{Problem formulation}

The volume of the patient's body is discretized into a three-dimensional grid of voxels. The voxels are assigned to different planning structures, e.g.\ organs or tumor tissue.
Let $\mathbf{P}$ denote the so-called dose matrix that maps a vector $x$ of
fluence intensities to the dose $d=\mathbf{P}x$ received by the individual voxels. Thus a planning structure $\mathcal{O}$ is essentially a subset of indices of the vector $d$.

To evaluate the dose received by different planning structures we used lower
and upper tail penalty functions as well as mean upper tail penalty
functions. The lower tail penalty function of the dose $d$ for a planning
structure $\mathcal{O}$ penalizes dose values below a given threshold $L$. It is given by

\begin{equation}  \label{eq:LowerTailPenalty}
f_{\mathcal{O},low}(d)\colonequals (1 / |\mathcal{O}|) \sum_{i\in \mathcal{O}}\left(
\max\{L-d_{i},0\}\right) ^{2}\quad L\in \mathbb{R}.
\end{equation}%
The upper tail penalty function is

\begin{equation}  \label{eq:UpperTailPenalty}
f_{\mathcal{O},up}(d)\colonequals (1 / |\mathcal{O}|) \sum_{i\in \mathcal{O}}\left( \max\{d_{i}-U,0\}\right)
^{2}\quad U\in \mathbb{R}
\end{equation}%
and is used to penalize dose values exceeding a given threshold $U$. The mean upper tail function is defined as

\begin{equation}  \label{eq:MeanUpperTailPenalty}
f_{\mathcal{O},mean}(d)\colonequals\left(\max\left\{ \left( (1 / |\mathcal{O}|) \sum_{i\in \mathcal{O}}d_{i}\right) -M, 0\right\}\right) ^{2}\quad
M\in \mathbb{R}.
\end{equation}%
This function is used to keep the mean dose values below a given threshold $M$, which is a less strict way of avoiding overdosage than using the upper tail penalty function.

Table \ref{tab:problemFormulation} shows our problem formulation for the
first head neck cases according to \cite{cco14}. The other cases were treated similarly, depending on the individual geometry of the tumor tissue. To gain an impression of the geometry of some of the planning structures involved see Figure \ref{fig:IMRTRSLOCT}. Note that the function values of the evaluation functions listed as constraints were required to be 0 which translates to the dose prescription of the corresponding planning structure being met by the dose distribution.

The weights for the evaluation functions were chosen empirically for each case individually such that the optimization would yield a clinically acceptable treatment plan.

Non-tumor Tissue 1 denotes the three-dimensional margin of 1 cm around the lymphatic drainage pathways (which provide the largest of the target volumes and contain both PTV 60 and PTV 70). Non-tumor Tissue 2 is the tissue that consists of everything but the lymphatic drainage pathways and Non-tumor Tissue 1.

These structures are introduced to reduce the radiation to non-tumor tissue, especially to areas which are not covered by any of the organs at risk. Additionally, it offers an incentive to the optimization algorithm to reduce the dose values of the organs at risk that lie outside of the lymphatic drainage pathways even below the prescribed maximum or maximum mean dose values.

In optimization level $\mu $ we minimize the objective function $\phi_{\mu}$ which is a weighted sum of the functions belonging to priority group $I_{\mu }$.

\begin{table}[tbp]
\caption{Evaluation functions of the structures under consideration gathered in different priority groups.}
\label{tab:problemFormulation}%
\begin{tabular}{lclc}
\hline
planning structure $\mathcal{O}$ & function & parameters & priority group \\ \hline
PTV 70 & lower tail penalty & $L$ = 66.5 & constraint \\
PTV 70 & upper tail penalty & $U$ = 80.5 & constraint \\
PTV 70 & mean upper tail penalty & $M$ = 73.5 & constraint \\
Myelon & upper tail penalty & $U$ = 45 & constraint \\ \hline
PTV 60 & lower tail penalty & $L$ = 57.5 & $I_1$ \\
Eye right & upper tail penalty & $U$ = 45 & $I_1$ \\
Eye left & upper tail penalty & $U$ = 45 & $I_1$ \\
Optic nerve right & upper tail penalty & $U$ = 50 & $I_1$ \\
Optic nerve left & upper tail penalty & $U$ = 50 & $I_1$ \\
Parotid left & mean upper tail penalty & $M$ = 26 & $I_1$ \\ \hline
Lymphatic drainage & lower tail penalty & $L$ = 47.5 & $I_2$ \\
Parotid right & mean upper tail penalty & $M$ = 26 & $I_2$ \\
Brain & upper tail penalty & $U$ = 60 & $I_2$ \\
Brainstem & upper tail penalty & $U$ = 52 & $I_2$ \\
Larynx & mean upper tail penalty & $M$ = 45 & $I_2$ \\
Esophagus & mean upper tail penalty & $M$ = 45 & $I_2$ \\
Oral cavity & mean upper tail penalty & $M$ = 40 & $I_2$ \\
Non-tumor Tissue 2 & mean upper tail penalty & $M$ = 5 & $I_2$ \\ \hline
Plexus right & upper tail penalty & $U$ = 63 & $I_3$ \\
Plexus left & upper tail penalty & $U$ = 63 & $I_3$ \\
Mandible & upper tail penalty & $U$ = 75 & $I_3$ \\
Inner ear right & upper tail penalty & $U$ = 45 & $I_3$ \\
Inner ear left & upper tail penalty & $U$ = 45 & $I_3$ \\
Lips & upper tail penalty & $U$ = 30 & $I_3$ \\
Non-tumor Tissue 1 & mean upper tail penalty & $M$ = 20 & $I_3$ \\\hline\hline
\end{tabular}%
\end{table}

\FloatBarrier

\subsubsection{Implementation}
\label{subsec:Implementation}
In this section we present pseudo codes of the algorithms we implemented to
solve the IMRT optimization problems. Algorithm \ref{Alg:1} implements the lexicographic optimization problem turned into a sequence of convex feasibility problems and Algorithm \ref{Alg:2} implements the superioriation methodology.

In the following paragraph we give some details about the submethods mentioned in the codes.\vspace\baselineskip

\texttt{findFeasibleSolution($t^{(\mu)}$)} tries to solve the CFP %
\eqref{eq:MCLOauxCFP} with $t_k^{(\mu)} = t^{(\mu)}$ within $n_{max}$ simultaneous
projections. If this is successful, it returns [$x$, true] with $
\phi_{\mu}(\mathbf{P}x)\leq t_k^{(\mu)}$ and $x$ feasible. If this is not successful,
the method returns [$x$, false] with $\phi_{\mu}(\mathbf{P}x)\leq
\tilde{t}^{(\mu)}$, where $\tilde{t}^{(\mu)}>t_k^{(\mu)}$ and $x$ feasible. $\tilde{t}^{(\mu)}$ is the smallest
upper bound for $\phi_{\mu}$ that the algorithm was able to find a feasible
solution for.

\texttt{reduce($\phi_{\mu}(\mathbf{P}x)$)} returns a new upper bound $t^{(\mu)}$ for
$\phi_{\mu}$ that is smaller than $\phi_{\mu}(\mathbf{P}x)$.

\texttt{addConstraint($\phi_{\mu}^{\ast}$)} adds $\phi_{\mu}(x) \leq
\phi_{\mu}^{\ast} + \delta_{\mu}$ to the set of constraints.

\texttt{supDirection($x$)} returns the direction of superiorization at $x$.
In our calculations we used superiorization with respect to the objective
function $\phi _{\mu +1}$ of the following optimization level.

\begin{algorithm}[H]
  \caption{Lexicographic level set scheme}\label{Alg:1}
  \begin{algorithmic}[1]

    \State $\mu$ = 1;
  \While{$\mu\leq M$}
    \If{$\mu==1$}
            \State $t^{(\mu)} = \infty$
        \Else
            \State $t^{(\mu)} = $ reduce($\phi_{\mu}(\textbf{P}x)$)
    \EndIf

        \State levelIsSolved = false
        \State counter = 0

        \While{$\neg$ levelIsSolved}

            \State [x, feasibleSolutionFound] = findFeasibleSolution($t^{(\mu)}$)
            \State counter++
            \If{feasibleSolutionFound}
                \If{$t^{(\mu)}==t_{min}$}
                \State levelIsSolved = true
                \State $\phi_{\mu}^{\ast} = \phi_{\mu}(\textbf{P}x)$
                \State addConstraint($\phi_{\mu}^{\ast}$)
                \State $\mu$++
                \Else
                \If{$\mu<M$ $\&\&$ counter $(\bmod K)$ == 0}
                \State x = superiorize()
                \EndIf
                \State $t^{(\mu)}$ = reduce($\phi_{\mu}(\textbf{P}x)$)
                \EndIf
            \Else
                \State levelIsSolved = true
                \State $\phi_{\mu}^{\ast} = \phi_{\mu}(\textbf{P}x)$
                \State addconstraint($\phi_{\mu}^{\ast}$)
                \State $\mu$++
            \EndIf
        \EndWhile
    \EndWhile
  \end{algorithmic}
\end{algorithm}

\begin{algorithm}[H]
  \caption{Superiorize}\label{Alg:2}
  \begin{algorithmic}[1]
    \State coefficientBigEnough = true
    \State base $\in (0, 1)$
    \State $\lambda$ = 0
    \While{$\lambda<\Lambda$ $\&\&$ coefficientBigEnough}
        \State direction = supDirection($x$)
        \State loop = true
        \State exponent = 1
        \While{loop $\&\&$ coefficientBigEnough}
            \State coefficient = base$^{\text{exponent}}$
            \State $x_{sup}$ = $x$ + coefficient $\cdot$ direction
            \If{$x_{sup}$ $\geq 0$ $\&\&$ $\phi_{\mu+1}(\textbf{P}x_{sup}) \leq \phi_{\mu+1}(\textbf{P}x)$}
                \State $x = x_{sup}$
                \State loop = false
                \State $\lambda$++
            \Else
                \If{coefficient $>$ minStepsize}
                    \State exponent++
                \Else
                    \State coefficientBigEnough = false
                \EndIf
            \EndIf
        \EndWhile
    \EndWhile
  \end{algorithmic}
\end{algorithm}

\section{Numerical illustration}

\label{Sec:Numerical illustration}

In the first part of this section we demonstrate the plausibility of using
the superiorized lexicographic optimization methodology on a
simple example with two optimization variables which is taken from Stanimirovi\'{c} \cite%
{Stanimirovic12}. In the second part we present the performance of our
method on four real IMRT head neck cases.

\subsection{2D example}

To demonstrate the mechanism of superiorization and how it can be useful for
lexicographic optimization we consider the following example.

\begin{example}
Solve the following multicriteria optimization problem

\begin{align}\label{ToyEx}
\min \text{ } &\Phi (x) =
\begin{pmatrix}
\phi _{1} \\
\phi _{2} \\
\phi _{3}
\end{pmatrix}
=
\begin{pmatrix}
-8x_{1}-12x_{2} \\
-14x_{1}-10x_{2} \\
-x_{1}-x_{2}
\end{pmatrix}
\notag \\
\text{such that  }
& 2x_{1}+x_{2}-150 \leq 0 \notag\\
& 2x_{1}+3x_{2}-300 \leq 0 \notag\\
& 4x_{1}+3x_{2}-360 \leq 0 \notag\\
& -x_{1}-2x_{2}+120 \leq 0 \notag\\
& -x_{1} \leq 0 \notag\\
& -x_{2} \leq 0
\end{align}
where the functions $\phi _{1},$ $\phi _{2}$ and $\phi _{3}$ are given in
lexicographic order.
\end{example}

The optimal solution $x^{\ast }=(x_{1}^{\ast },x_{2}^{\ast })$ is $(30,80)$
with \newline
$\Phi (x^{\ast }) = (\phi _{1}(x^{\ast }),\phi _{2}(x^{\ast }), \phi
_{3}(x^{\ast })) = (-1200,-1220,-110)$.

We stopped the optimization process at the iterate $x^{k}$ if $\Vert \Phi
(x^{k})-\Phi (x^{\ast })\Vert \leq 10^{-2}$. In the following, points are
given with precision $10^{-2}$. In this example calculating the gradients
and projections is computationally very cheap, but for IMRT cases gradient-
and dose-evaluations consume most of the computation time, so we compare the
number of projections and gradient evaluations.\vspace\baselineskip

\begin{figure}[tbp]
    \centering
    \subfigure[unsuperiorized LO]
    {\includegraphics[width=0.4\textwidth, height =0.4\textwidth]{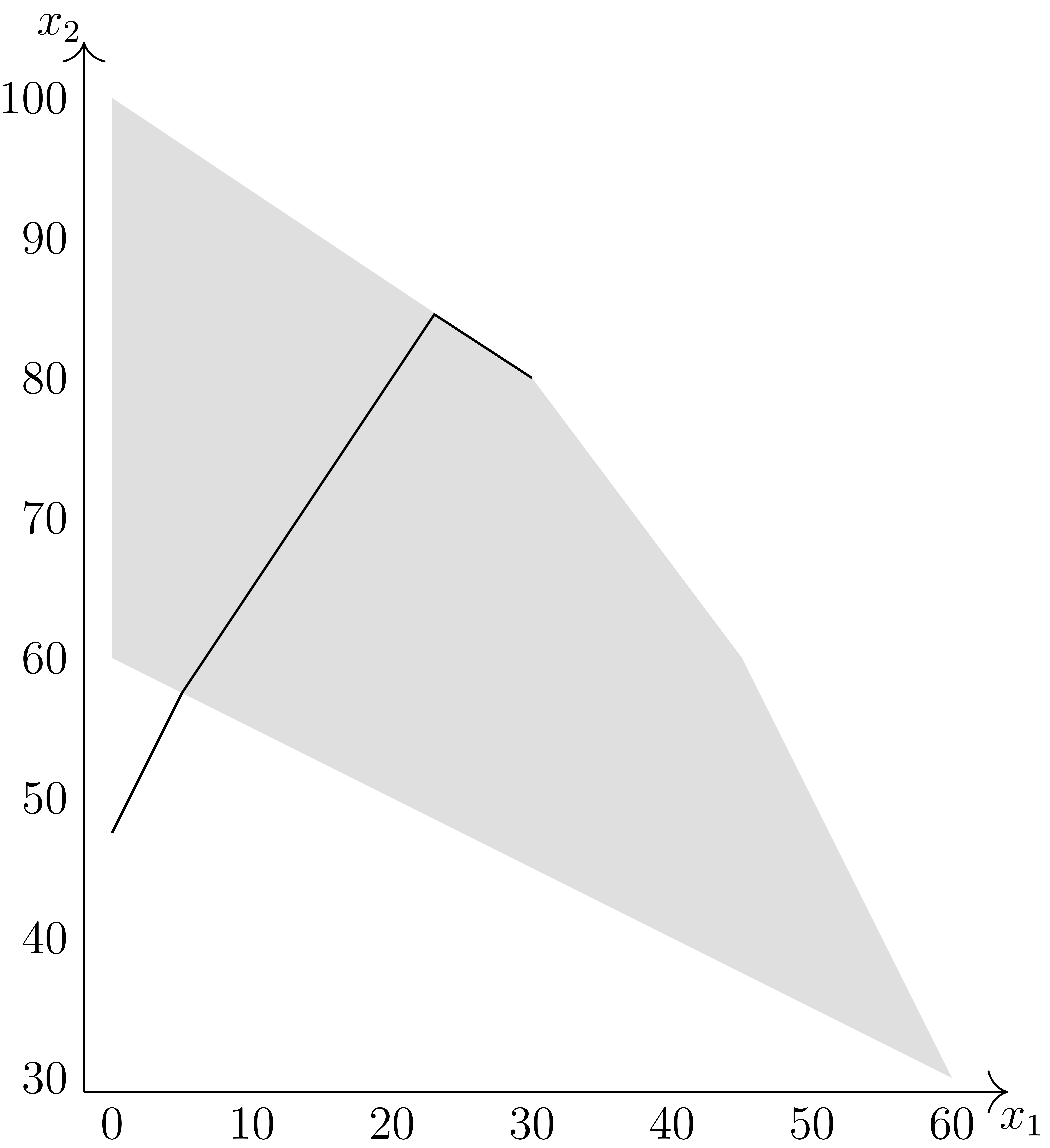}
    \label{fig:ToyClassicLO}}
    \subfigure[superiorized LO]
    {\includegraphics[width=0.4\textwidth, height =0.4\textwidth]{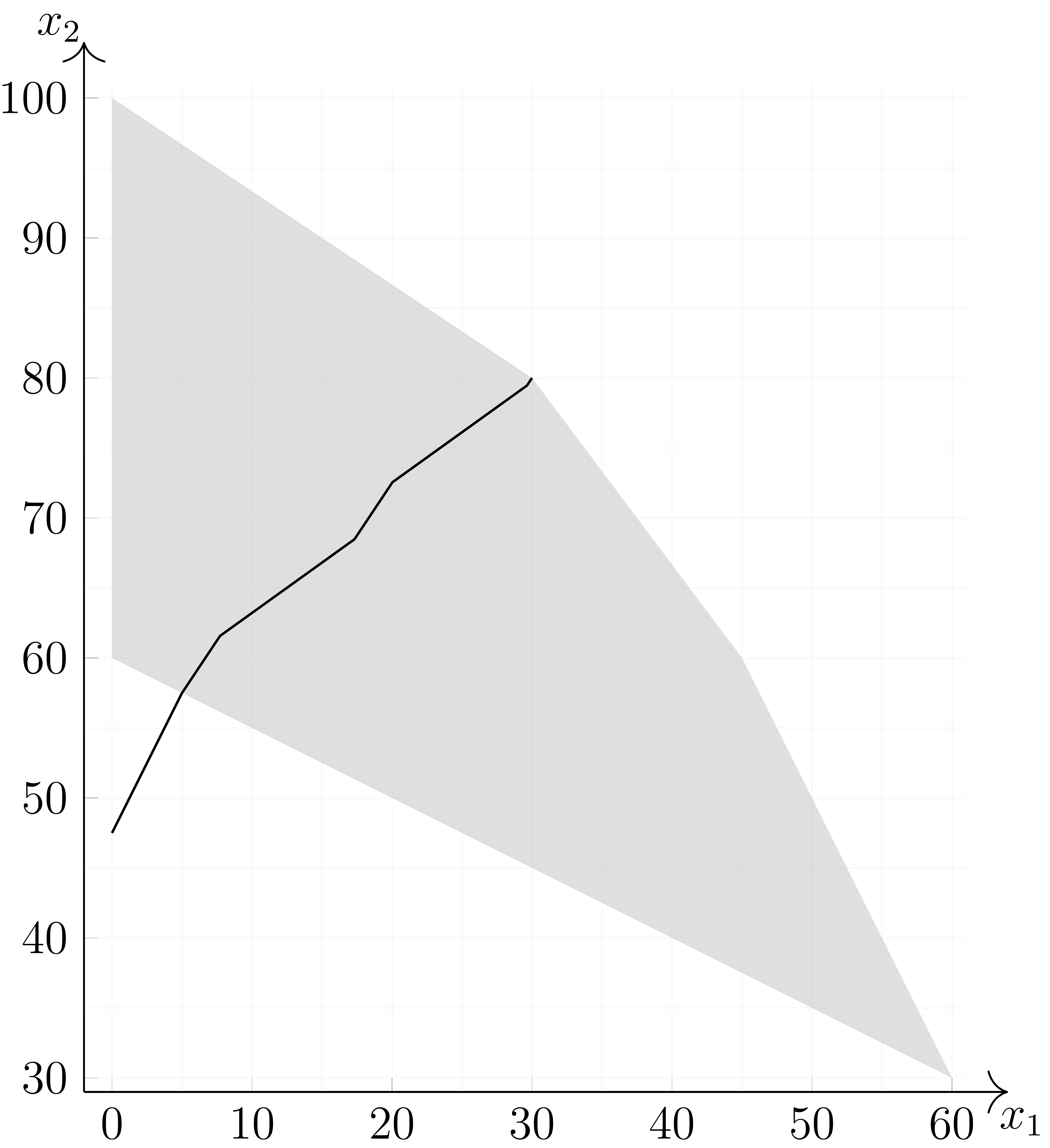}
    \label{fig:ToySupLO}}
    \caption{Trajectories for unsuperiorized and superiorized LO. The set of feasible solutions is shown in light gray.}
    \label{fig:ToyExTrajectoryBoth}
\end{figure}

Figure \ref{fig:ToyExTrajectoryBoth} shows the feasible region and the
trajectories of the iterates of the classical and the superiorized LO for
the starting point $x^{0}=(0,47.5)$.

Both methods first seek feasibility and find the point $(5,57.5)$. In
optimization level 1, the classical level set scheme minimizes with respect
to $\phi_{1}$ only and ends up at the solution $x^{(1)}=(23.08,84.625)$ of
the first optimization level with $\phi_{1}^{\ast }=\phi _{1}(x^{(1)})=-1200
$ using 327 projections and 327 gradient calculations. Additionally, we need
another 260 projections and 518 gradient evaluations to find out that we
cannot improve $\phi_{1}$ any further while staying feasible. In
optimization level 2, we want to preserve $\phi_{1}(x^{k})=\phi_{1}^{\ast } $
for all subsequent iterates $x^{k}$ and minimize $\phi_{2}$ under this
additional condition. We arrive at the solution $x^{(2)}=(30,80)$ of optimization level 2  after a total of 4743 projections and 9107 gradient
calculations. We cannot improve any further with respect to $\phi_{3}$ and
have found the optimal solution.\vspace\baselineskip

The superiorized level set scheme minimizes with respect to $\phi_{1}$ in
optimization level 1 but uses superiorization with respect to $\phi_{2}$. It
arrives at the solution $x^{(1)} = (30, 80)$ of optimization level 1 after
108 projections and 144 gradient calculations, 34 of which were needed for
the superiorization. Using the superiorized LO we get already sufficiently
close to the optimum in optimization level 1.

\subsection{IMRT}

From the formulae \eqref{eq:LowerTailPenalty}--\eqref{eq:MeanUpperTailPenalty} of the utilized functions it is clear that
the minimal value we can hope to achieve for any of the $\phi _{\mu }$ is 0.
Thus we considered an optimization level $\mu $ to be solved, if $\phi _{\mu
}(\mathbf{P}x)\leq t_{min}\colonequals 10^{-8}$ given that $x$ is feasible.\vspace\baselineskip

What takes most of the computational time with IMRT cases are
multiplications of the dose matrix $\mathbf{P}$ with the fluence vector $x$
which occur when we evaluate the dose or calculate the objective function
gradient for a given $x$. The number of these multiplications is what we
will use as units to measure progress, i.e., reduction of objective function
values, of the classical and the superiorized LO.

In Algorithm \ref{Alg:1} and \ref{Alg:2} in Subsection \ref{subsec:Implementation} we find parameters $K$ and $\Lambda$ which determine the behavior of the method. $K$ is the number of CFPs that have to be successfully solved before superiorization is applied and $\Lambda$ is the maximum number of superiorization steps that are taken. If the superiorization step size is smaller than a predefined threshold or the maximum number of steps is reached, the superiorization stops and we solve the next CFP.

In Tables \ref{tab:IMRTfastResults1} and \ref{tab:IMRTfastResults2} we first
present results we achieved by trying different parameter sets $(K,\Lambda )$
and then picking for each case individually the set that resulted in the
fastest convergence of the algorithm measured in the number of
multiplications. We refer to the sets of parameters picked according to this
strategy as $(\hat{K},\hat{\Lambda})$. All values presented in the tables in this section are rounded to $10^{-4}$.

The results presented in Table \ref{tab:IMRTfastResults1} show that fast SLO results in equal or lower objective values than unsuperiorized LO within less multiplications. In Table \ref{tab:IMRTfastResults2} we see that fast SLO produces solutions of the single optimization levels with potentially lower objective values of the subsequent objective functions than unsuperiorized LO. Fast SLO thus offers a better starting point to the single optimization levels than unsuperiorized LO which reduces the overall number of multiplications.\vspace\baselineskip

$\phi _{\mu }^{\ast }=\phi _{\mu }(x^{(\mu )})$ denotes the value of
$\phi _{\mu }$ at the solution $x^{(\mu )}$ of optimization level $\mu $
obtained using the classical LO, $\hat{\phi}_{\mu }^{\ast }=\phi _{\mu }(%
\hat{x}^{(\mu )})$ denotes the value of $\phi _{\mu }$ at the solution $\hat{%
x}^{(\mu )}$ of optimization level $\mu $ obtained using the superiorized LO
with $(\hat{K},\hat{\Lambda})$.

$N_{j}^{k}$ is the number of multiplications needed by the classical LO to
solve optimization level $j+1$ to level $k$. $N\colonequals N_{0}^{M}=%
\sum_{j=0}^{M-1}N_{j}^{j+1}$ is the total number of multiplications.
Analogously, $\hat{N}_{j}^{k}$ denotes the number of multiplications needed
by the superiorized LO with $(\hat{K},\hat{\Lambda})$ to solve optimization
level $j+1$ to level $k$ and $\hat{N}\colonequals \sum_{j=0}^{M-1}\hat{N}_{j}^{j+1}$.

For our calculations we chose $\delta _{\mu }=0.1\phi _{\mu }^{\ast }$
for the classical LO and $\delta _{\mu }=0.1\hat{\phi}_{\mu }^{\ast }$
for the superiorized LO.

\begin{table}[H]
\caption{\small Optimal objective value and total multiplication number ratios for $(\hat{%
K}, \hat{\Lambda})$. Fast SLO achieves lower or equal objective values than unsuperiorized LO within fewer multiplications. The notation $1^{\lozenge }$ means that we solved the optimization
level, meaning that the objective function value has dropped below $10^{-8}$.}
\label{tab:IMRTfastResults1}
\begin{center}
\begin{tabular}{lcccc}
\hline
\rule{0pt}{15pt} & $\hat{\phi}_1^{\ast}/\phi_1^{\ast}$ & $\hat{\phi}%
_2^{\ast}/\phi_2^{\ast}$ & $\hat{\phi}_3^{\ast}/\phi_3^{\ast}$ & $\hat{N}$/$%
N $ \\ \hline
case 1 & $1^{\lozenge}$ & 0.9598 & 0.9923 &  0.8998\\
case 2 & $1^{\lozenge}$ & 0.9743 & 0.9960 & 0.6612 \\
case 3 & $1^{\lozenge}$ & 0.9996 & 1.0000 & 0.9210 \\
case 4 & $1^{\lozenge}$ & 0.9708 & 0.9941 & 0.8959 \\ \hline\hline
\end{tabular}%
\end{center}
\end{table}

\begin{table}[H]
\caption{\small Multiplication numbers per optimization level and objective values
of starting points for $(\hat{K}, \hat{\Lambda})$. Fast SLO produces solutions of the single optimization levels with potentially lower objective values of the subsequent objective functions than unsuperiorized LO.}
\label{tab:IMRTfastResults2}%
\begin{tabular}{lccccc}
\hline
\rule{0pt}{15pt} & $\hat{N}_0^1/N_0^1$ & $\hat{N}_1^2/N_1^2$ & $\hat{N}%
_2^3/N_2^3$ & $\phi_2(\hat{x}^{(1)})/\phi_2(x^{(1)})$ & $\phi_3(\hat{x}%
^{(2)})/\phi_3(x^{(2)})$ \\ \hline
case 1 & 1.6521 & 0.5669 & 0.9895 & 0.3741 & 0.9856 \\
case 2 & 1.4066 & 0.3588 & 1.0065 & 0.3682 & 1.0008 \\
case 3 & 1.0403 & 0.8594 & 1.0000 & 0.5625 & 1.0075 \\
case 4 & 0.3818 & 1.0536 & 0.8591 & 0.8196 & 0.9939 \\ \hline\hline
\end{tabular}%
\end{table}

\FloatBarrier

It is not clear a priori how to choose values for $K$ and $\Lambda$ for a
given problem. Therefore we aimed to find a robust set of parameters that
yields good results for all IMRT cases and could possibly be applied for
other IMRT head neck cases as well. We refer to this set of parameters as $(%
\check{K}, \check{\Lambda})$. The results for the robust parameter set are
presented in Table \ref{tab:IMRTrobustResults1} and \ref%
{tab:IMRTrobustResults2}.

The results in Table \ref{tab:IMRTrobustResults1} show that robust SLO achieves comparable objective values to unsuperiorized LO within fewer multiplications. In Table \ref{tab:IMRTrobustResults2} we see that robust SLO, too, produces solutions of the single optimization levels with potentially lower objective values of the subsequent objective functions than unsuperiorized LO. Like fast SLO, robust SLO also offers a better starting point to the single optimization levels than unsuperiorized LO which again reduces the overall number of multiplications.\vspace\baselineskip

Following the previous notation $\check{\phi}_{\mu }^{\ast }=\phi_{\mu }(%
\check{x}^{(\mu )})$ is the value of $\phi _{\mu }$ at the solution $\check{x%
}^{(\mu )}$ of optimization level $\mu $ obtained using the superiorized LO
with $(\check{K}, \check{\Lambda})$.

$\check{N}_{j}^{k}$ denotes the number of multiplications needed by the
superiorized LO with $(\check{K}, \check{\Lambda})$ to solve optimization
level $j+1$ to level $k$ and $\check{N} \colonequals \sum_{j=0}^{M-1}\check{N}%
_{j}^{j+1}$.

\begin{table}[H]
\caption{\small Optimal objective value and total multiplication number ratios for $(%
\check{K}, \check{\Lambda})$. Robust SLO achieves comparable objective values to unsuperiorized LO within fewer multiplications. The notation $1^{\lozenge }$ means that we solved the optimization
level, meaning that the objective function value has dropped below $10^{-8}$.}
\label{tab:IMRTrobustResults1}
\begin{center}
\begin{tabular}{lcccc}
\hline
\rule{0pt}{15pt} & $\check{\phi}_1^{\ast}/\phi_1^{\ast}$ & $\check{\phi}%
_2^{\ast}/\phi_2^{\ast}$ & $\check{\phi}_3^{\ast}/\phi_3^{\ast}$ & $\check{N}
$/$N$ \\ \hline
case 1 & $1^{\lozenge}$ & 1.0046 & 1.0002 & 0.9177 \\
case 2 & $1^{\lozenge}$ & 1.0020 & 1.0043 & 0.8877 \\
case 3 & $1^{\lozenge}$ & 0.9829 & 0.9969 & 0.9541 \\
case 4 & $1^{\lozenge}$ & 1.0245 & 1.0071 & 0.8859 \\ \hline\hline
\end{tabular}%
\end{center}
\end{table}

\begin{table}[H]
\caption{\small Multiplication numbers per optimization level and objective values
of starting points for $(\check{K}, \check{\Lambda})$. Robust SLO produces solutions of the single optimization levels with potentially lower objective values of the subsequent objective functions than unsuperiorized LO.}
\label{tab:IMRTrobustResults2}%
\begin{tabular}{lccccc}
\hline
\rule{0pt}{15pt} & $\check{N}_0^1/N_0^1$ & $\check{N}_1^2/N_1^2$ & $\check{N}%
_2^3/N_2^3$ & $\phi_2(\check{x}^{(1)})/\phi_2(x^{(1)})$ & $\phi_3(\check{x}%
^{(2)})/\phi_3(x^{(2)})$ \\ \hline
case 1 & 1.1825 & 0.6972 & 0.9849 & 0.3916 & 1.0008 \\
case 2 & 0.7447 & 0.8409 & 0.9999 & 0.5194 & 1.0058 \\
case 3 & 1.1320 & 0.9029 & 0.9999 & 0.5530 & 1.0007 \\
case 4 & 1.2453 & 0.8575 & 0.8664 & 0.4959 & 1.0103 \\ \hline\hline
\end{tabular}%
\end{table}

\FloatBarrier

\begin{figure}[tbp]
\begin{minipage}{\textwidth}
\centering
\includegraphics[height = 0.45\textheight]{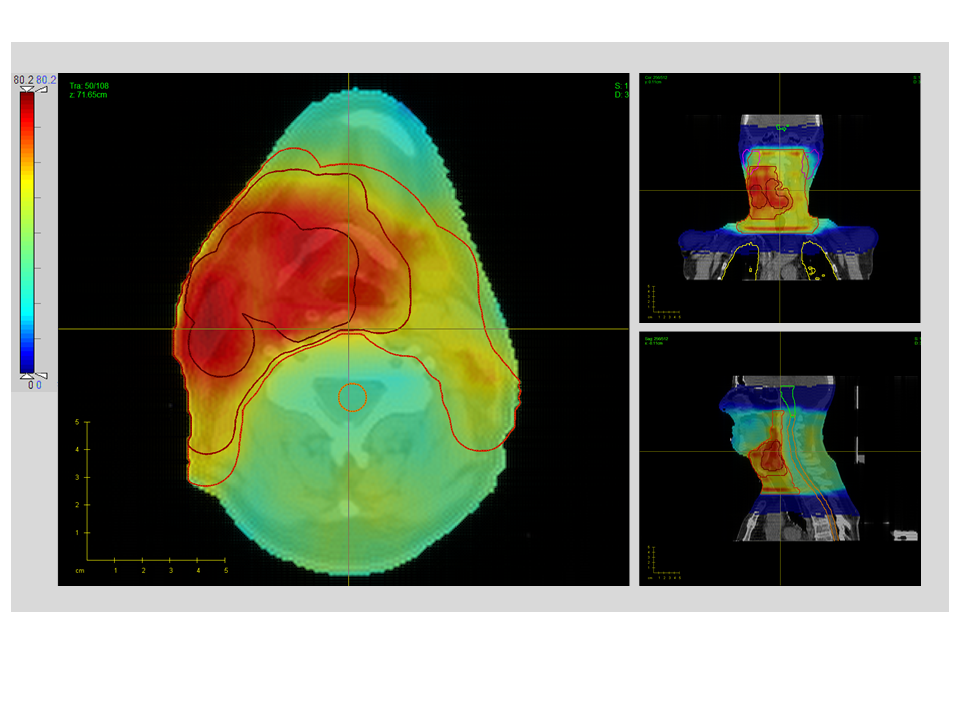}
\caption{CT image and colorwash of dose resulting from robustly superiorized LO. The contours of the target volumes PTV 70, PTV 60 and the lymphatic drainage pathways (LDP) are shown in dark, medium and bright red. The contours of the left and right parotid are marked in dark and bright magenta. The lungs are contoured in yellow and the myelon is shown in orange. Finally, the brainstem is contoured in bright light green.}
\label{fig:IMRTRSLOCT}

\includegraphics[height = 0.35\textheight, width = \textwidth]{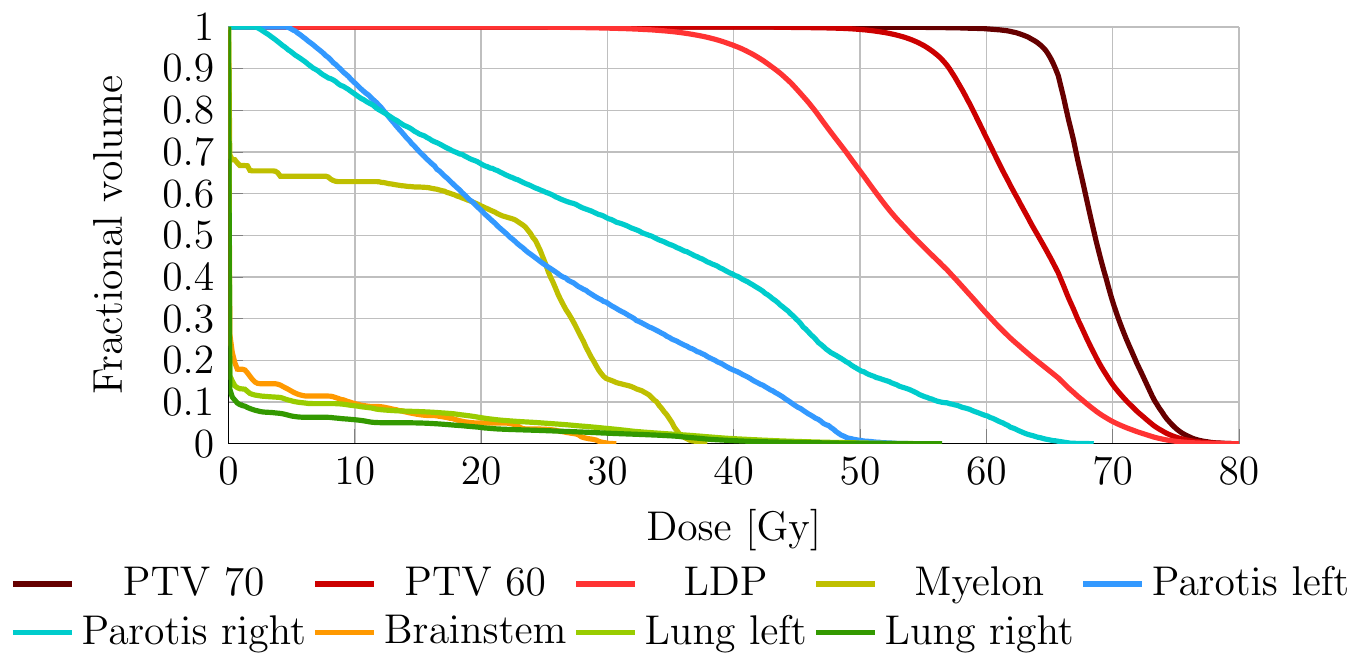}
\caption{DVHs resulting from solution dose of robustly superiorized LO.}
\label{fig:IMRTRSLODVH}
\end{minipage}
\end{figure}

In Figure \ref{fig:IMRTRSLOCT} we present one slice of the CT of case 1 with the dose distribution resulting from the solution of superiorized LO with the robust parameter set $(\check{K},\check{\Lambda})$. The corresponding DVH is shown in Figure \ref{fig:IMRTRSLODVH}.

Both the images and the dose matrix were created using CERR \cite{cerr03}. Note that the dose calculation via the dose matrix, which is used in the optimization algorithm, provides only an approximation to the dose values calculated by CERR which are depicted in Figure \ref{fig:IMRTRSLOCT} and \ref{fig:IMRTRSLODVH}. Thus the DVH curves show e.g. that dose values below 66.5 Gy occur in the PTV70 even though in the optimization the mathematical constraint ensuring that this does not happen is not violated.

\FloatBarrier

\begin{figure}[tbp]
\centering
\subfigure[Values of $\phi_1$ for robust SLO, fast SLO and unsuperiorized LO.] {\includegraphics[width = \textwidth, height = 0.23\textheight] {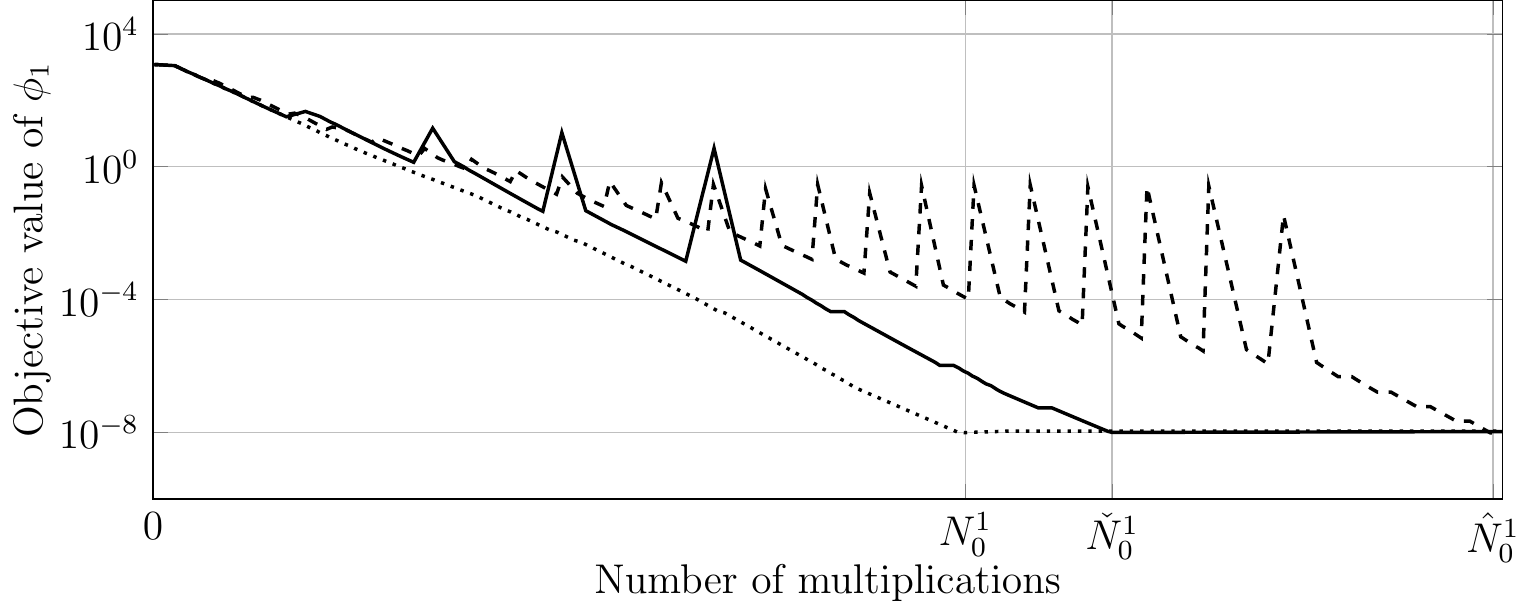}\label{fig:IMRTobjfunctionsA}}
\\
\subfigure[Values of $\phi_2$ for robust SLO, fast SLO and unsuperiorized LO.] {\includegraphics[width = \textwidth, height = 0.23\textheight] {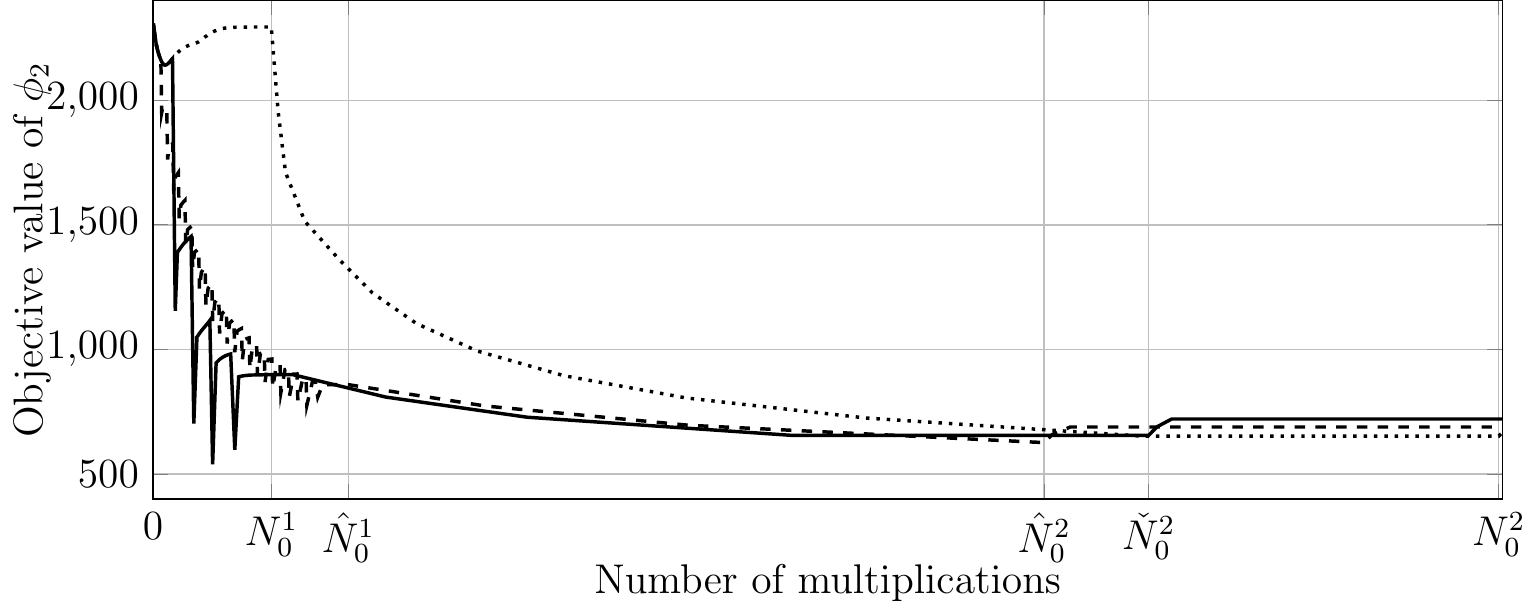}\label{fig:IMRTobjfunctionsB}}
\\
\subfigure[Values of $\phi_3$ for robust SLO, fast SLO and unsuperiorized LO.] {\includegraphics[width = \textwidth, height = 0.23\textheight] {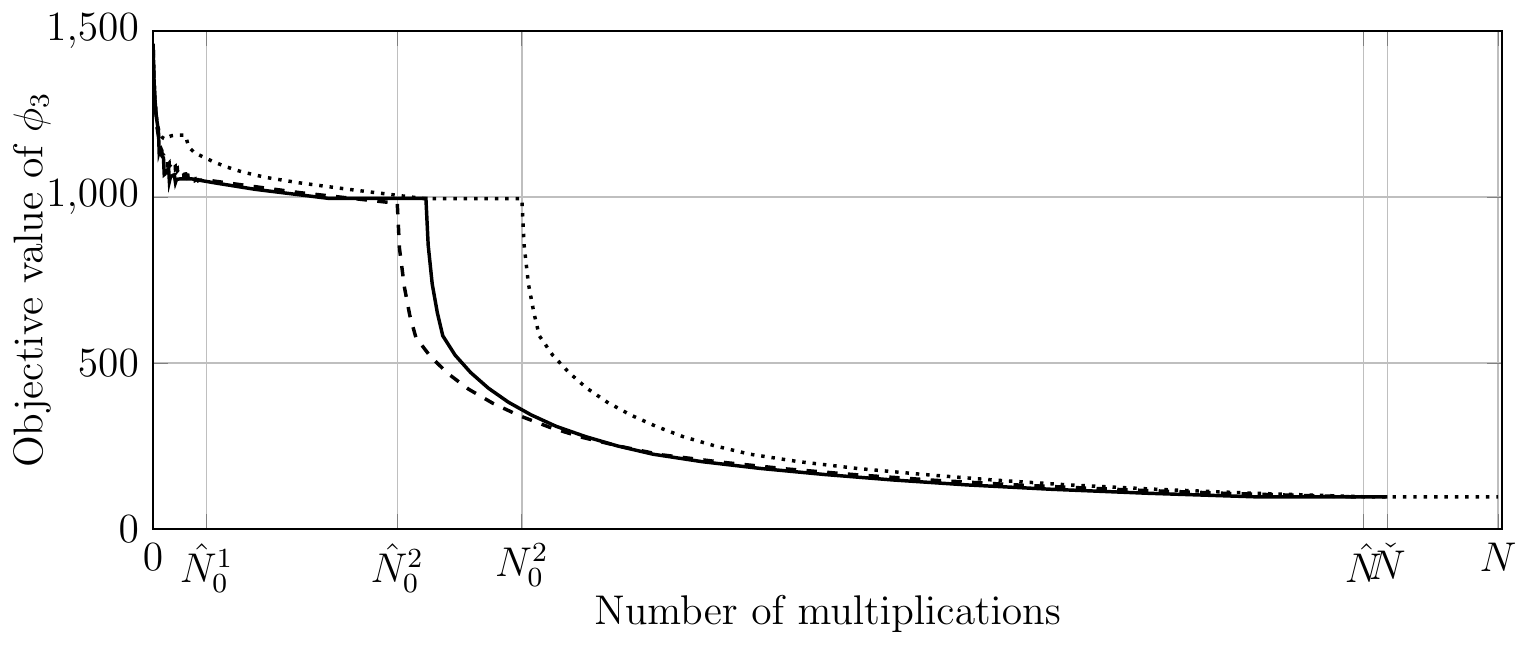}\label{fig:IMRTobjfunctionsC}}\\
\subfigure
{\includegraphics[height = 0.03\textheight]
{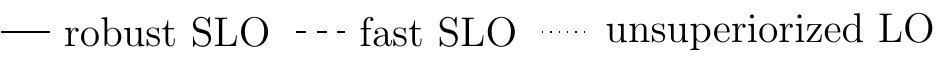}}\\
\caption{Superiorization perturbs the optimization of $\phi_1$. Superiorization applied in level 1 offers a starting point with significantly lower value of $\phi_2$ and speeds up the optimization of $\phi_2$. No superiorization is applied in level 2 and SLO and unsuperiorized LO perform similarly. Due to the speedup in level 2 the overall number of multiplications is lower for SLO.}
\label{fig:IMRTobjfunctions}
\end{figure}

In Figure \ref{fig:IMRTobjfunctions} we present the objective values of $
\phi _{1},\phi _{2}$ and $\phi _{3}$ plotted against the number of
dosematrix-vector-multiplications for case 1. To show the fundamental behavior of the superiorization method we intentionally did not give any actual numbers in the plots in Figure \ref{fig:IMRTobjfunctions}.

In Figure \ref{fig:IMRTobjfunctionsA} we can see how $\phi _{2}$ increases
while we minimize $\phi _{1}$ in level 1 using classical LO. We see how the
superiorization with respect to $\phi _{2}$ disturbs the minimization of $
\phi _{1}$ and leads to $N_{0}^{1}$ being smaller than $\hat{N}_{0}^{1}$.
Figure \ref{fig:IMRTobjfunctionsB} shows that $\hat{x}^{(1)}$ is a better
starting point for the minimization of $\phi _{2}$ in level 2 than $x^{(1)}$ because its objective value of $\phi_2$ is significantly lower. This results in $\hat{N}_{1}^{2}$ being notably smaller than $N_{1}^{2}$.

Due to the correlation of the clinical goals there are relatively few
auxiliary CFPs to be solved in optimization level 2 for all cases. When
superiorization is applied (both with optimal and robust $(K,\Lambda )$) we
even get so close to the minimum of $\phi_{2}$ at the end of level 1 that
for these parameter pairs the CFP in level 2 turns inconsistent within less than $\hat{K}$ or $\check{K}$ reductions of $t_k^{(2)}$. Therefore no
superiorization is applied in level 2 and the objective value of $\phi_3$ is about the same for fast SLO, robust SLO and unsuperiorized LO.

As we can see from Figure \ref{fig:IMRTobjfunctionsC} superiorization with
respect to $\phi _{2}$ during level 1 can in some cases -- depending on the
correlation of the objective functions -- reduce the value of $\phi_{3}$ and $
\hat{N}_{2}^{3}$.

Similar behavior can be observed for the superiorized LO that uses the robust parameter set $(\check{K}, \check{\Lambda})$.

\section{Discussion}
\label{Sec:Discussion}

We used LO as an intuitive approach to lower the complexity of decision making in IMRT and to include a priori knowledge about the priorities of the decision maker. We combined this MCO technique with the superiorization methodology to speed up the optimization in the next optimization level.
This is not the only application for superiorization. One could for example choose as superiorization function $\psi$ a weighted sum of different evaluation functions than we did, or one could choose $\psi$ as a function that is not included in the lexicographic optimization model at all.\vspace\baselineskip
In this work we identified the superiorization parameters manually. This process takes of course too much time to be of any practical use.
To be able to find robust parameters that work for many optimization problems in an automated way, one should make sure to identify problems with similar tradeoffs among the objective functions and the superiorization function(s).\vspace\baselineskip

An issue with using projection methods to solve CFPs is that in general we do not have a criterion that tells us whether a CFP is
inconsistent. Therefore we defined a maximum number $n_{max}$
of attempts, i.e.\ simultaneous projections, to solve it. If the CFP cannot
be solved within $n_{max}$ simultaneous projections, we assume it to be
inconsistent. This means that finding out that a CFP actually is inconsistent takes $n_{max}$ projections, which drastically increases the total number of projections.

On the other hand, low values of $n_{max}$ might not suffice to reach the actual minimum of an objective function under the given constraints, but only an approximation thereof.

We decided to stop solving an optimization level $\mu$ after encountering the first inconsistent auxilary CFP \eqref{eq:MCLOauxCFP} and continue with level $\mu+1$, accepting that we might miss the minimal value of $\phi_{\mu}$ accessible with the chosen value of $n_{max}$ by a certain gap (which is bounded by the reduction strategy of $t^{(\mu)}$, in our case by 10 $\%$).

With this strategy, the computationally expensive process of finding out that a CFP is inconsistent occurs (at most) once per optimization level.

For this reason it is a big advantage to know the objective values of the optimal solution a priori as in Example \ref{ToyEx}. In general we do not have this
information but as in our calculations for the IMRT cases, sometimes it is
possible to exploit a priori information about the objective functions to
define a stopping criterion that is computationally cheap to evaluate.\vspace\baselineskip

A question arising in the context of optimization problems is how algorithm 1 performs compared to other solvers. In the following, we will refer to algorithm 1 as (superiorized) projection solver.
We implemented the projection solver in Matlab. To keep computation time for algebraic operations comparable we chose to compare our solutions to the solutions of Matlab's fmincon function using the interior point solver.

We did our calculations in Matlab 2015b on an Intel Core i7-5600U processor with 2.6 GHz and 8GB RAM.

Matlab's interior point solver took 16904 to 30690 seconds to solve the lexicographic optimization problem for our four cases. Within a maximum of 1000 iterations per optimization level it yielded solutions which were slightly infeasible (maximal constraint violation of $\mathcal{O}(10^{-8})$ -  $\mathcal{O}(10^{-6})$). Due to the nature of the evaluation functions this infeasibility had no major effect on the dose values. These solutions were still acceptable from a medicinal point of view.

For a detailed comparison of computation time see Table \ref{tab:MatlabSolverComparisonCompTime}.

\begin{table}[H]
    \caption{\small Computation time in seconds used by Matlab's interior point solver (IP), the unsuperiorized projection solver, the projection solver using robust parameters $(\check{K}, \check{\Lambda})$ (RobustSProj), and the projection solver using fast parameters $(\hat{K}, \hat{\Lambda})$ (FastSProj).}
    \label{tab:MatlabSolverComparisonCompTime}
    \begin{center}
        \begin{tabular}{lcccc}
            \hline
            \rule{0pt}{13pt}& IP & Proj & RobustSProj & FastSProj \\\hline
            case 1 & 30689.8 & 55.3 & 50.8 & 49.8\\
            case 2 & 17649.8 & 14.4 & 12.8 & 9.5\\
            case 3 & 26596.2 & 13.9 & 13.6 & 12.9\\
            case 4 & 16903.9 & 14.8 & 13.0 & 13.5\\ \hline\hline
        \end{tabular}%
    \end{center}
\end{table}

All iterates produced by the projection solver are feasible by construction of the algorithm. This enables us to decide how much computational effort we want to put into the calculations - of course at the expense of mathematical solution quality.

For our calcuations we chose the maximum number $n_{max} = 10^3$ of simultaneous projections, to solve a given CFP for th projection solver. When comparing the resulting solutions to the ones produced by Matlab's solver (see Tables \ref{tab:MatlabSolverComparisonObjValues} and \ref{tab:ProjSolverComparisonObjValues}) it becomes apparent that they are suboptimal.
Solutions with equal or smaller objective function values can be found by the projection solver, too, if $n_{max}$ is raised.

\begin{table}[H]
    \caption{\small Objective funtion values (rounded to precision $10^0$ if $\geq$ 1) at the solution $x^{\ast}$ found by Matlab's interior point solver}
    \label{tab:MatlabSolverComparisonObjValues}
    \begin{center}
        \begin{tabular}{lcccc}
            \hline
            \rule{0pt}{13pt}& $\phi_1(x^{\ast})$ & $\phi_2(x^{\ast})$ & $\phi_3(x^{\ast})$ & maximal constraint \\
            & & & & violation \\\hline
            \rule{0pt}{15pt}case 1 & 9.1529$\cdot 10^{-8}$ & 425 & 145 & 9.1419$\cdot10^{-6}$\\
            case 2 & 3.6419$\cdot10^{-7}$ & 363 & 795 & 7.2779$\cdot10^{-6}$\\
            case 3 & 0.3834 & 172 & 371 & 9.8338$\cdot10^{-7}$\\
            case 4 & 0.4778 & 66 & 382 & 9.6511$\cdot10^{-8}$\\ \hline\hline
        \end{tabular}%
    \end{center}
\end{table}

\begin{table}[H]
    \caption{\small Objective funtion values (rounded to precision $10^0$ if $\geq$ 1) at the solution $x^{\ast}$ found by the superiorized projection solver with robust superiorization parameters $(\check{K}, \check{\Lambda})$ and $n_{max} = 10^3$.}
    \label{tab:ProjSolverComparisonObjValues}
    \begin{center}
        \begin{tabular}{lcccc}
            \hline
            \rule{0pt}{13pt}& $\phi_1(x^{\ast})$ & $\phi_2(x^{\ast})$ & $\phi_3(x^{\ast})$ & maximal constraint  \\
            & & & & violation  \\\hline
            \rule{0pt}{15pt}case 1 & $\leq 1.1\cdot 10^{-8}$ & 720 & 97 & $\leq 10^{-8}$\\
            case 2 & $\leq 1.1\cdot 10^{-8}$ & 502 & 1725 & $\leq 10^{-8}$\\
            case 3 & $\leq 1.1\cdot 10^{-8}$ & 1477 & 1918 & $\leq 10^{-8}$\\
            case 4 & $\leq 1.1\cdot 10^{-8}$ & 1318 & 851 & $\leq 10^{-8}$\\ \hline\hline
        \end{tabular}%
    \end{center}
\end{table}

\begin{figure}[H]
    \centering
    \includegraphics[width = \textwidth]{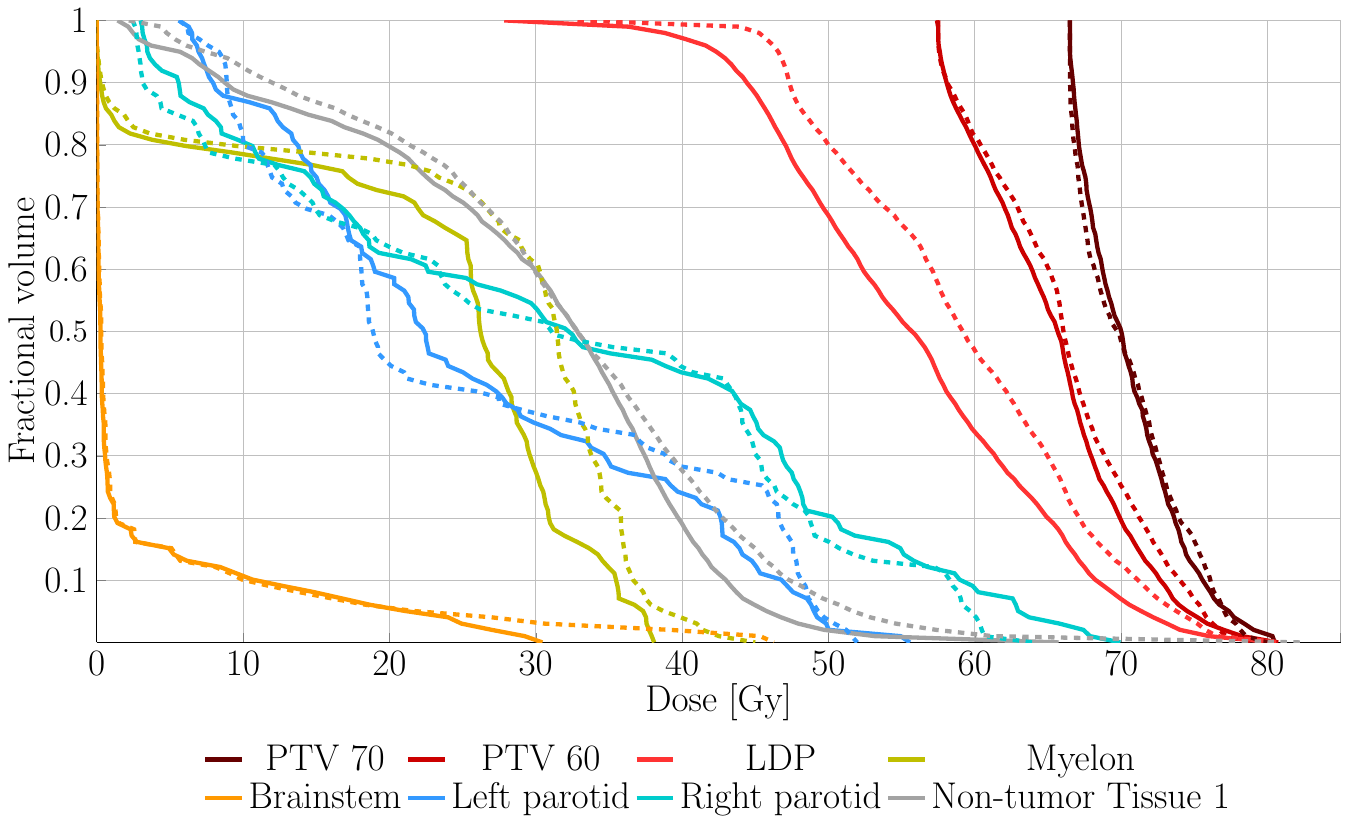}
    \caption{DVH comparison of solution doses of superiorized projection solver using superiorization parameters $(\check{K}, \check{\Lambda})$ and Matlab's interior point solver (dashed) for case 1.}
    \label{fig:IMRTSolverComparison}
\end{figure}

The DVH curves of case 1 in Figure \ref{fig:IMRTSolverComparison} show that the additional 30639 seconds of computation time the Matlab solver takes compared to RobustSProj result in an improvement of 3 Gy in the average dose value of the LDP (whose dose evaluation function is contained in $I_2$) at the cost of 2 Gy increase in the average dose value of non-tumor tissue 1 (whose dose evaluation function is contained in $I_3$ and thus considered less important) .

For case 1, the Matlab solver yields a solution with has a lower evaluation function value corresponding to the LDP than the projection solver. Because in case 1 the LDP has a relatively big weight $w_i$ in \eqref{eq:Phi_mu}, this difference is amplified.
We observed similar behavior with the other IMRT cases.

Concerning the question of performance compared to other solvers, we observe that the projection solver with $n_{max} = 10^3$ offers a fast way to reach a solution that is in many aspects already as good as the solution produced by Matlab's solver. If one is willing to invest additional hours of computation time, however, Matlab's solver offers solutions with lower objective function values.

\section{Conclusion}

\label{Sec:Conclusion}

We demonstrated the concept of superiorization in combination with an MCO method on a simple example as well as on four IMRT head neck cases. We showed that for LO we can speed up the optimization process by providing a better starting point for the next level of optimization. We found a common set of superiorization parameters which yielded robust results for all four cases.

\section{Acknowledgments}

We thank the Department of Radiation Oncology (Head: Prof. Dr. C. Belka) of the Ludwig-Maximilians University (LMU), Munich, Germany, for providing the clinical data for our calculations and visualizations.

\end{document}